\documentclass[11pt, psamsfonts]{amsart}
\usepackage{amssymb}
\usepackage{array,calc,url,tabularx, epic,tikz}
\usetikzlibrary{arrows,positioning,decorations.pathmorphing,
  decorations.markings}

\def\ignore#1{\relax}

\def\g{\mathfrak g}
\def\h{\mathfrak h}

\def\sp{{\mathfrak sp}}
\def\sl{{\mathfrak sl}}

\def\Z{{\mathbb Z}}

\def\C{{\mathbb C}}

\def\la{\lambda}
\def\La{\Lambda}

\def\S{\mathcal S}
\def\B{\mathcal B}

\def\Ca{\mathcal C}
\def\D{\mathcal D}
\def\Na{\mathcal N}

\def\E{\mathcal E}

\def\M{\mathcal M}
\def\B{{\mathcal B}}

\def\ignore#1{\relax}

\def\1{{\bf 1}}

\def\End{{\rm End}}
\def\Hom{{\rm Hom}}

\def\HnN{H_n^{(N)}}
\def\varep{\varepsilon}

\calclayout

\renewenvironment{proof}[1][\proofname]{\par
  \normalfont
  \topsep6\p@\@plus6\p@ \trivlist
  \item[\hskip\labelsep\bfseries
    #1\@addpunct{.}]\ignorespaces
}{%
  \qed\endtrivlist
}

\def\th@plain{%
  \let\thmhead\thmhead@plain \let\swappedhead\swappedhead@plain
  \thm@preskip.5\baselineskip\@plus.2\baselineskip
                                    \@minus.2\baselineskip
  \thm@postskip\thm@preskip
  \itshape
\renewcommand{\labelenumi}{{(\alph{enumi})\quad}}
                        \renewcommand{\labelenumii}{{(\roman{enumii})\ }}
}
\def\th@definition{%
  \let\thmhead\thmhead@plain \let\swappedhead\swappedhead@plain
  \thm@preskip.5\baselineskip\@plus.2\baselineskip
                                    \@minus.2\baselineskip
  \thm@postskip\thm@preskip
  \upshape
}
\def\th@remark{%
  \thm@headfont{\itshape}
  \let\thmhead\thmhead@plain \let\swappedhead\swappedhead@plain
  \thm@preskip.5\baselineskip\@plus.2\baselineskip
                                    \@minus.2\baselineskip
  \thm@postskip\thm@preskip
  \upshape
}

{\theoremstyle{plain}
\newtheorem{theorem}{Theorem}[section]
}

{\theoremstyle{plain}
\newtheorem{proposition}[theorem]{Proposition}
}

{\theoremstyle{plain}
\newtheorem{corollary}[theorem]{Corollary}
}

{\theoremstyle{plain}
\newtheorem{lemma}[theorem]{Lemma}
}

{\theoremstyle{plain}

}

{\theoremstyle{definition}
\newtheorem{definition}[theorem]{Definition}
}

{\theoremstyle{definition}

}

{\theoremstyle{remark}
\newtheorem{remark}[theorem]{Remark}
}

{\theoremstyle{remark}

}

\numberwithin{equation}{section}

\renewcommand{\labelenumi}{{ \theenumi.}}
\renewcommand{\labelenumii}{{(\alph{enumii})}}

\def\eth{e_{(3)}}

\def\la{\lambda}
\def\al{\alpha}

\def\vo{v_o}
\def\mnl{m_{n,\la}}

\def\choose #1 #2{\begin{pmatrix}#1\\#2\end{pmatrix}}

\def\Np{[N]_{+}}
\def\Nm{[N]_{-}}
\def\CnN{C_n^{(N)}}

\ignore{
\input BoxedEPS
\SetRokickiEPSFSpecial
\HideDisplacementBoxes
\SetEPSFDirectory{./graphics/} }

\ignore{
\input BoxedEPS
\SetTexturesEPSFSpecial
\HideDisplacementBoxes
\SetEPSFDirectory{:graphics:}
}

\begin{document}

\title[Module categories]
{On module categories related to $Sp(N-1)\subset Sl(N)$}

\author{ Hans Wenzl}

\address{Department of Mathematics\\ University of California\\ San Diego,
California}

\email{hwenzl@ucsd.edu}

\begin{abstract}
Let $V=\C^N$ with $N$ odd. We construct a $q$-deformation of $\End_{Sp(N-1)}(V^{\otimes n})$  which contains
$\End_{U_q\sl_N}(V^{\otimes n})$. It is a quotient of an abstract two-variable algebra which is defined by adding
one more generator to the generators of the Hecke algebras $H_n$. These results suggest the existence of module categories of  $Rep(U_q\sl_N)$
which may not come from already known coideal subalgebras of $U_q\sl_N$. We moreover indicate how this can be used to construct module categories
of the associated fusion tensor categories as well as subfactors, along the lines of previous work for
inclusions $Sp(N)\subset SL(N)$.
\end{abstract}
\maketitle
\centerline{IN MEMORIAM OF VAUGHAN JONES WITH GRATITUDE}
\vskip .5cm

The problem of classifying module categories of a given tensor category arises in different contexts
such as conformal field theory and the study of subfactors. A lot of progress has been made for module categories
of fusion categories coming from quantum groups (or Wess-Zumino-Witten models) in recent work by Edie-Michell
(see \cite{Cain}), building on work of Jones,  Ocneanu, Gannon, Schopieray, Evans and Pugh  and others. But in many cases, a detailed
description of those module categories such as fusion rules, algebras etc is still not available.

It is well-known that  $Rep(H)$ is a module category of $Rep(G)$ for an inclusion of groups $H\subset G$.
The basic idea here is to find a subgroup $H$ of a Lie group $G$ for which one can find analogs of the module
$Rep(H)$ in fusion categories related to $G$. This was successfully carried out in  \cite{WqBrauer} and \cite{Wfusion}
for inclusions  $Sp(N)\subset SL(N)$, $N$ even, and $O(N)\subset SL(N)\times \Z_2$
for arbitrary $N>1$. It allowed a detailed description of these module categories in terms of the well-known combinatorics
of these groups (see Section \ref{conclusions}).   The current paper roughly contains the analogous results of the paper \cite{WqBrauer} 
for the inclusions $Sp(N-1)\subset SL(N)$ for $N$ odd.

Here is our approach in more detail, formulated for the group $G=SL(N)$ for simplicity.
We denote the quotient tensor category of tilting modules of the quantum group $U_q\sl_N$ modulo
negligible modules for $q$ a root of unity by $\overline{Rep}(U_q\sl_N)$. It is often also denoted by $SU(N)_k$
for $q^2$ a primitive $(N+k)^{th}$ root of unity.
 We study the following questions: Find a subgroup $H\subset SL(N)$
for which we can find

(a) a $q$-deformation of $Rep(H)$ which is a module category of $Rep(U_q\sl_N)$,

(b) a quotient of said $q$-deformation for $q$ a root of unity which is a module category of  

$\overline{Rep}(U_q\sl_N)$,

(c) a subfactor corresponding to the module category in (b) if it is unitarizable.

\noindent  Before describing the results in this paper in more detail,
we would like to make a few general remarks about this approach. It is known that Question (a) can be solved if $H$ is
the group of fixed points under a period 2 automorphism via coideal subalgebras, see work of Letzter \cite{Le} and of Noumi and Sugitani \cite{NS}.
It follows from work in \cite{Mo} that  the examples in \cite{Wfusion} indeed correspond to special cases of the work in \cite{Le} and \cite{NS}.
However, it is also clear that not all module categories constructed via coideal algebras allow solutions of Questions (b) and (c).
This can be seen e.g. for $SL(3)$, where all module categories of the corresponding fusion categories are known due to work of
Gannon and Evans and Pugh, see \cite{EP}.

As stated in the title, we consider the inclusion of $Sp(N-1)\subset SL(N)$ for $N$ odd. One of the main results of this  paper
is the construction of a sequence of  two-parameter algebras $C_n=C_n^{(N)}(q)$ which contain the Hecke algebras $H_n(q)$
of type $A_{n-1}$ as subalgebras. For $q\neq \pm 1$ it suffices to add one more generator $e$, corresponding to the projection
onto the trivial $Sp(N-1)$ submodule $\1\subset V$, where $V=\C^N$ is the vector representation of $SL(N)$. In order to get the correct algebra in the classical limit $q\to \pm 1$, 
we also need additional generators to obtain nontrivial morphisms 
between the two copies of the trivial representation in $V^{\otimes 2}$.
 We give  a presentation of these algebras via generators and relations and an explicit basis.
Another important result is the proof of the existence of an extension of the Markov trace for $H_n$ to the algebras
$C_n$ which satisfies a generalized Markov condition.

We now give a more detailed description of the results of this paper which will also explain how these algebras can be used to
construct module categories.
We fix notations in the first section and prove a number of combinatorial and algebraic results concerning  
$\End_{Sp(N-1)}(V^{\otimes n})$. We derive relations for a $q$-deformation
of  $\End_{Sp(N-1)}(V^{\otimes n})$ which contains the Hecke algebra $H_n$ as a subalgebra in the second section.
These relations
are essentially forced by the fusion respectively restriction rules for $Sp(N-1)\subset SL(N)$ and a generalized Markov
condition (see Condition \ref{Markovprop} or the discussion here for Section 4).  In particular, we show that
only two solutions are possible for fixed $q\neq \pm 1$ and $N>1$ odd
(see Remark \ref{modrestriction}). Hence we have two possible choices $C_{n,\pm}$
of extensions of $H_n$ subject to our conditions. As they are closely related, see Remark \ref{relationsremark},
we will often just use the notation $C_n$ for either of these cases. We also define a version of $C_n$ depending on two
parameters $q$ and $p$, which specializes to the original version for $p=q^N$. It is then shown that for $N>2n$
the dimension of these algebras is at most $\dim \End_{Sp(N-1)}(V^{\otimes n})$ and  explicit spanning sets are determined.
We define representations of $C_n^{(N)}$ into $\End(V^{\otimes n})$ in the third section. It is shown that their images
in the classical limit $q=1$ concide with $\End_{Sp(N-1)}(V^{\otimes n})$. We conclude from this that the given spanning
sets in the previous section are actually bases. 
Section 4 contains a proof that we can extend the Markov trace $tr$ of the Hecke algebra $H_n$ to a trace on $C_n$
satisfying the generalized Markov condition $tr(cg_n)=tr(c)tr(g_n)$ for all $c\in C_n$.
In the last section, we first briefly describe how this paper has been influenced by the work of Vaughan Jones, even though
this will be pretty obvious to experts anyways. We then indicate how our algebras can be used to construct module categories
and subfactors. To do so, we consider the quotients $\overline{H}_n$ and $\overline{C}_n$ of the algebras $H_n$ and $C_n$
modulo the annihilator ideals of the Markov trace $tr$.  The objects of the category and of the module category are given by idempotents
in $\overline{H}_n$ and $\overline{C}_n$ respectively, and the module action comes from the natural inclusion map
$C_n\otimes H_m\to C_{n+m}$.
To illustrate this, we restate results from \cite {WqBrauer} and \cite{Wfusion} in the language of module
categories which was not used there. This shows, in particular, that Problems (a)-(c) have been solved for the inclusions
considered in those papers. We finally discuss how the approach in this paper can be used to give detailed descriptions 
for a large class of module categories.

$Acknowledgements:$ I would like to thank Gail Letzter and Cain Edie-Michell for useful conversations.

\section{Fusion rules for the embedding of $Sp(N-1)\subset SL(N)$, $N$ odd}

\subsection{Fusion rules} Let $V=\C^N$ with $N$ odd. We fix a symplectic bilinear form $(\ ,\ )'$
on $V$ with 1-dimensional kernel spanned by the nonzero vector $\vo$ and a complement $V'$ of $\vo$ on which the 
form $(\ ,\ )'$ is nondegenerate. This defines an  embedding $Sp(N-1)\subset SL(N)$ such that $V$
 decomposes into the direct sum $V'\oplus \C\vo$ as
an $Sp(N-1)$ module. 
Recall that the finite-dimensional simple representations of $Sp(N-1)$ are labelled by Young diagrams with $\leq (N-1)/2$
rows, with the $(N-1)$-dimensional simple representation $V'$ labeled by the Young diagram with one box.
 If $V_\la$ is a simple representation of $Sp(N-1)$ labeled by the Young diagram $\la$, we have 
\begin{equation}\label{tensorsymp}
V_\la\otimes V'\cong \bigoplus_{\mu\leftrightarrow\la} V_\mu,
\end{equation}
where $\mu$ ranges over all Young diagrams which can be obtained from $\la$ by either adding or removing a box.
One deduces from this the tensor product rule
\begin{equation}\label{tensor}
V_\la\otimes V\cong V_\la\oplus \bigoplus_{\mu\leftrightarrow\la} V_\mu,
\end{equation}
with $\mu$ as in \ref{tensorsymp}. 

\subsection{Bratteli diagrams and path bases}\label{Bratteli}
 The inclusions of the algebras 
$$... \subset \End_{Sp(N-1)}(V^{\otimes n})
 \subset \End_{Sp(N-1)}(V^{\otimes n+1})\subset \ ... $$
 are conveniently described by a Bratteli diagram. It follows from the 
tensor product rules \ref{tensor}  that its vertices at level $n$ are labeled  by Young diagrams $\la$ with $|\la|\leq n$,
where $|\la|$ denotes the number of boxes in the Young diagram. A diagram $\la$ at level $n$ is connected with a diagram
$\mu$ at level $n+1$ if $\la$ differs from $\mu$ by at most one box. The multiplicity of $V_\la$ in $V^{\otimes n}$ is then
given by the number of paths of length $n$ from level 0 to level $n$. As we only have multiplicities 0 or 1, we obtain a basis
of $\Hom(V_\la, V^{\otimes n})$ labeled by the paths of length $n$ which end in $\la$ for any irreducible $Sp(N-1)$ module
$V_\la$ labeled by $\la$. This basis is unique up to rescaling by nonzero scalars.  Below is the part of the Bratteli diagram containing level 1
and level 2.

\vskip 0in 
\begin{picture}(400,120)(0,0)  
\put(140,50){\circle*{4}}
\put(180,50){\circle*{4}} 
\put(220,50){\circle*{4}} 
\put(220,50){\line(-1,1){40}}
\put(260,50){\circle*{4}} 
\put(260,50){\line(-2,1){80}}
\put(140,50){\line(0,1){40}}
\put(140,50){\line(1,1){40}}
\put(140,90){\circle*{4}}
\put(180,90){\circle*{4}}
\put(180,50){\line(0,1){40}}
\put(180,50){\line(-1,1){40}}
\put(137,35){$\emptyset$}
\put(175,35){[1]}
\put(214,35){[2]}
\put(253,35){[1,1]}
\put(137,97){$\emptyset$}
\put(175,97){[1]}
\put(145,0){Inclusion diagram}
\put(310,50){$\End_{Sp(N-1)}(V^{\otimes 2})$}
\put(310,90){$\End_{Sp(N-1)}(V)$}
\end{picture}
\vskip 0.5cm

\subsection{Multiplicities for large $N$}
The tensor product rule \ref{tensor} allows us to calculate the multiplicity $\mnl$ of the simple module $V_\la$ in $V^{\otimes n}$. 
If $N>2n$, there are no restrictions on Young diagrams, and we can give closed formulas for the multiplicities $\mnl$.
To do so, we define integers $h_r$ inductively by $h_0=1$, $h_1=1$ and
\begin{equation}\label{hr}
h_{r+1}= h_r + rh_{r-1}.
\end{equation}
We denote  by $d_\la$ the dimension of the simple $S_n$ module labeled by the Young diagram $\la$, where 
the number of boxes $|\la|$ of $\la$ is equal to $n$. There exists a well-known explicit formula for it in terms of
the hook lengths of $\la$. We will need the following well-known identities
\begin{equation}\label{dlambdas}
\sum_{\mu<\la} d_\mu\ =\ d_\la\quad {\rm and}\quad \sum_{\nu>\la} d_\nu = (n+1)d_\la,
\end{equation}
where $|\la|=n$, and $\mu$ and $\nu$ range over all Young diagrams which can be obtained by removing a box
from $\la$ (for $\mu$) or adding a box to $\la$ (for $\nu$).

\begin{proposition}\label{multiplicities} (a) The multiplicity $\mnl$ of the simple module $V_\la$ in $V^{\otimes n}$ for $N>2n$
 is equal to
$h_{n-|\la|}\binom{n}{|\la|}d_\la$.

(b) We have the identity
$$\sum_{|\la|=n-r}\mnl d_\la\ =\ h_r\frac{n!}{r!}$$

(c) If $N>2n$, we have $\dim \End_{Sp(N-1)}(V^{\otimes n})  =  \sum_{r=0}^n h_r^2\frac{n!}{r!}\binom{n}{r}$.
\end{proposition}

$Proof.$ We will prove (a) by induction on $n$, with the claim easily checked for $n=1$. It follows from \ref{tensor} that
$V_\la\subset V_\mu\otimes V$ if and only if $\mu=\la$ or $\mu$ is obtained by removing or adding a box from/to $\la$. In each of theses 
cases $V_\la$ appears with multiplicity 1. It follows that
$$\mnl\ =\ m_{n-1,\la}\ +\ \sum_{\mu<\la} m_{n-1,\mu}\ +\ \sum_{\nu>\la} m_{n-1,\nu}.$$
Using the induction assumption, the identities \ref{hr}, \ref{dlambdas} and the identity
$\binom{n-1}{k+1}(k+1)=\binom{n-1}{k} (n-1-k)$, we obtain
\begin{align} \mnl\ &= h_{n-|\la|-1}\binom{n-1}{|\la|}d_\la \ +\ \sum_{\mu<\la} h_{n-|\la|}\binom{n-1}{|\la|-1}d_\mu
\ +\ \sum_{\nu>\la} h_{n-|\la|-2}\binom{n-1}{|\la|+1}d_\nu\notag\\
&= h_{n-|\la|-1}\binom{n-1}{|\la|}d_\la\ +\ h_{n-|\la|}\binom{n-1}{|\la|-1}d_\la\ +\ h_{n-|\la|-2}\binom{n-1}{|\la|+1}d_\la(|\la|+1)\notag\\
&=\ d_\la[h_{n-|\la|-1}\binom{n-1}{|\la|}\ +\ h_{n-|\la|}\binom{n-1}{|\la|-1}\ +\ h_{n-|\la|-2}\binom{n-1}{|\la|}(n-1-|\la|)\notag\\
&=\ d_\la[\binom{n-1}{|\la|}(h_{n-|\la|-1}+h_{n-|\la|-2}(n-1-|\la|))\ +\ h_{n-|\la|}\binom{n-1}{|\la|-1}\notag\\
&=\ d_\la h_{n-|\la|}( \binom{n-1}{|\la|} + \binom{n-1}{|\la|-1}).
\notag\end{align}
This proves part (a). Part (b) follows from this and the identity $\sum_{|\la|=n-r} d_\la^2=(n-r)!$. Part(c) similarly follows from this, part (a) and
$ \dim \End_{Sp(N-1)}(V^{\otimes n}) = \sum_{|\la|\leq n} \mnl^2$.

\subsection{ Description of $\End_{Sp(N-1)}(V^{\otimes n})$}
We denote by $E\in \End(V)$ the projection onto the trivial representation $\C \vo$ of $Sp(N-1)$ with kernel
$V'$, and by $U$  the antisymmetrization map
$$U:\ v\otimes w\ \mapsto\ v\wedge w = v\otimes w - w\otimes v, \quad v,w\in V.$$
The elements $U_i\in \End(V^{\otimes n})$ are defined for $1\leq i<n$ by
\begin{equation}\label{Hecketensor}
U_i\ =\ 1\otimes 1\otimes\ ...\ \otimes U\ \otimes \ ...\ \otimes 1,
\end{equation}
where $U$ acts on the $i$-th and $(i+1)$-st factors in $V^{\otimes n}$.
Observe that the flip $G: v\otimes w\in V^{\otimes 2}\mapsto w\otimes v$ is related to $U$ by
the simple formula $G=1-U$.
We now extend the bilinear form $(\ ,\ )'$ to a non-degenerate bilinear form $(\ ,\ )$ on $V$ by defining
\begin{equation}\label{nondeg}
(\vo,\vo)=1\quad {\rm and}\quad (\vo,v')=0=(v',\vo)\quad {\rm for \ all}\ v'\in V'.
\end{equation}
If $(v_i')$ and $(w_i')$ are dual bases of $V'$ with respect to $(\ ,\ )'$, we obtain the canonical 
vector $\vo\otimes\vo + \sum_i v_i'\otimes w_i'\in V^{\otimes 2}$ for the form $(\ ,\ )$. We define
the element $F\in \End(V^{\otimes 2})$ by
$$F(v\otimes w)\ =\ (v,w)\ (\vo\otimes\vo + \sum_i v_i'\otimes w_i').$$
The elements $F_i$ and $G_i$ in $\End(V^{\otimes n})$ are defined in the same way
as $U_i$  in \ref{Hecketensor}. It is well-known that the element $F$ can be used to calculate the trace
of an element $A\in \End(V)$ by
$$Tr(A)F=F(A\otimes 1)F.$$

Let us now decompose $V^{\otimes n}$ into the direct sum of three $Sp(N-1)$ submodules as follows.
Write
$V^{\otimes n}=\bigoplus_\la m_\la V_\la$,
where $m_\la$ is the multiplicity of the simple module $V_\la$ in $V^{\otimes n}$.
For given $n$, we call the diagram $\la$ an {\it old / recent / new diagram} in $V^{\otimes n}$
 if $V_\la$ has appeared for the first time 
in $V^{\otimes m}$ with $m\leq n-2$ / $m=n-1$ / $m=n$. 
We define 
$$V^{\otimes n}_{old}\ =\ \bigoplus_{\la\ old} m_\la V_\la,$$
and $V^{\otimes n}_{rec}$ and $V^{\otimes n}_{new}$ accordingly. 
Then it follows from the tensor product rule \ref{tensor} that
$\la$ is an old / recent/ new diagram if $|\la|\leq n-2$/$|\la| = n-1$ / $|\la|=n$. 

\begin{theorem}\label{basicconstruction} (see e.g. \cite{Wexc}, Proposition 4.10)
Let $V$ be a finite-dimensional self-dual $G$-module. Then $\End_G(V^{\otimes n}_{old})$
is given by a Jones basic construction for $\End_G(V^{\otimes n-2})\subset \End_G(V^{\otimes n-1})$.
Moreover, it coincides with the two-sided ideal in $\End_G(V^{\otimes n})$ generated by $F_{n-1}$, which 
is spanned by elements of the form $aF_{n-1}b$, with $a,b\in \End_G(V^{\otimes n-1}\otimes 1)$.
\end{theorem}

\begin{theorem}\label{fundamental}
The algebra $\End_{Sp(N-1)}(V^{\otimes n})$ is generated by the symmetric group $S_n$, acting via permutations
of the factors of $V^{\otimes n}$, the element $E\otimes 1_{n-1}$ and the element $F_1$.
\end{theorem}

$Proof.$ The proof goes by induction on $n$, with $n=1$ obviously true. The claim follows for $\End(V^{\otimes n})_{old}$
from Theorem \ref{basicconstruction} and the induction assumption. If $V_\la\subset V^{\otimes n}_{new}$,
i.e. $|\la|=n$, it follows from the tensor product rules that
$V_\la\subset (V')^{\otimes n}$. In this case, the claim follows from Brauer's classical result. It also implies that
$W_\la=\Hom(V_\la, V^{\otimes n})$ is an irreducible $S_n$ module.
To deal with the remaining cases, observe  that the elements $G_i$, $1\leq i<n$
and $T=(1-2E)\otimes 1_{n-1}$ satisfy the relations of the Weyl group $W(B_n)$ of type $B_n$. Hence the quotient modulo the
ideal generated by $F_{n-1}$ is also a quotient of the group algebra of $W(B_n)$. We will finish the proof in
the next section after a brief review of the representation theory of $W(B_n)$.

\subsection{Weyl group of type $B_n$} The Weyl group $W(B_n)$ of type $B_n$ 
can be defined via generators $t$ and $s_i$, $1\leq i< n$ 
and relations such that the $s_i$ are simple reflections of the symmetric group (e.g. we can take $s_i=(i,i+1)$) and such that
$t$ commutes with $s_i$ for $i>1$ and satisfies $s_1ts_1t=ts_1ts_1$. It is well-known that it is isomorphic to the
semidirect product $(\Z/2)^n\rtimes S_n$, with $t$ corresponding to $(1,0,...,0)\in (\Z/2)^n$, and with $S_n$ permuting the coordinates of 
elements of $(\Z/2)^n$. The irreducible representations of $W$ are labeled by pairs of Young diagrams $(\la,\mu)$ with 
$|\la|+|\mu|=n$. They can be constructed as follows (see e.g. \cite{Ser}, Chapter 8 for details):

Let $\phi$ be a character of  $(\Z/2)^n$ such that
$\phi(\varep_i)=-1$ for $i\leq r$ and $\phi(\varep_i)=1$  for $i>r$; here $\varep_i$ is the $i$-th unit vector in $(\Z/2)^n$.
The centralizer of $\phi$ consists of all $g\in W(B_n)$ such $\phi(gxg^{-1})=\phi(x)$ for all $x\in (\Z/2)^n$.
It is easy to see that for our choice of $\phi$ the centralizer is equal to $(\Z/2)^n\rtimes (S_r\times S_{n-r})$.
Let $W_\la$ and $W_\mu$ be irreducible representations of $S_r$ and $S_{n-r}$. Then $W_\la\otimes W_\mu$ 
becomes an irreducible representation of $(\Z/2)^n\rtimes (S_r\times S_{n-r})$, where the action of $x\in (\Z/2)^n$ is
given by the scalar $\phi(x)$.
It can then be shown that inducing this representation up to $W(B_n)$ yields an irreducible representation 
of $W(B_n)$ of dimension $\binom{n}{r}d_\la d_\mu$. 

\medskip

{\it Conclusion of  proof of Theorem \ref{fundamental}} Let $|\la|=n-1$. Then $E\otimes E\otimes 1_{n-2}$, and hence also $F_1$
acts as 0 on $\Hom(V_\la, V^{\otimes n})$, with $V_\la$ an irreducible $Sp(N-1)$-module.
Hence we can view  $\Hom(V_\la, V^{\otimes n})$ as a $W(B_n)$-module
on which $E\otimes 1_{n-1}$ acts nontrivially; indeed, the module $\C\vo\otimes V_\la\cong V_\la$ is in the image of 
$E\otimes 1_{n-1}$. As $\C v_0\otimes V_\la\subset \C v_o\otimes (V')^{\otimes n-1}$, it also follows that $E_i$
acts as 0 on $\C\vo\otimes V_\la$ for $i>1$. Hence the action of $(\Z/2)^n$ on $\C v_o\otimes V_\la$ is given by the functional
$\phi: x\in (\Z/2)^n\mapsto (-1)^{x_1}$.
 We obtain that  $\Hom(V_\la, V^{\otimes n})$ contains an irreducible $W(B_n)$-module labeled by 
$([1],\la)$. By the previous discussion it has dimension $n\ d_\la= w_{n,\la}=\dim \Hom(V_\la, V^{\otimes n}) $. Hence $W(B_n)$ and therefore also
$C_n$ acts irreducibly on   $\Hom(V_\la, V^{\otimes n})$.

\section{A $q$-deformation of $\End_{Sp(N-1)}(V^{\otimes n})$}

The goal of this paper is to study $q$ deformations of $Rep(Sp(N-1))$  at the categorical level which are compatible with 
the deformation of $Rep(Sl(N))$ to $U_q\sl_N$ and with the embedding  $Sp(N-1)\subset Sl(N)$. As we shall see,
this leads to a structure different from $Rep(U_q\sp_{N-1})$.

\subsection{Hecke algebras} The Hecke algebra $H_n=H_n(q)$ of type $A_{n-1}$ is defined via generators $g_i$, $1\leq i<n$ and
relations
$$g_ig_{i+1}g_i\ =\  g_{i+1}g_ig_{i+1},\quad 1\leq i<n-1,$$
together with $g_i^2=(q-q^{-1})g_i+1$ and $g_ig_j=g_jg_i$ for $|i-j|\neq 1$.
Defining 
$$u_i=q1-g_i,$$
the relations above translate to $u_i^2=(q+q^{-1})u$, $u_iu_j=u_ju_i$ for $|i-j|\neq 1$ and
\begin{equation}\label{Hecke3}
u_iu_{i+1}u_i-u_i\ =\ u_{i+1}u_iu_{i+1}-u_{i+1}.
\end{equation}

We replace $Rep(Sl(N))$ by the representation category of the Drinfeld-Jimbo quantum group
$U_q\sl_N$. It was shown in \cite{Ji} that the generalization of the map $U$ of the previous subsection can then be defined with respect to a basis $\{ v_i\}$ for $V$ by
$U(v_i\otimes v_i)=0$ and by defining its restriction to the ordered basis vectors $v_i\otimes v_j$ and $v_j\otimes v_i$, $i<j$ by the matrix
\begin{equation}\label{Heckerep}
 \ U\ =\ 
\left[\begin{matrix} q^{-1}&-1\\-1&q \end{matrix}\right].
\end{equation}
The elements $U_i\in \End(V^{\otimes n})$ are defined for $1\leq i<n$ as in \ref{Hecketensor}.

\subsection{Dimensions and traces}\label{tracesec} We review some basics about dimension functions for representations
of quantum groups, see \cite{Turaev} and  \cite{Ks}, XIV.4 for more details. Let $W$ be an $U_q\sl_N$-module
with dual module $W^*$. Then there exist canonical morphisms
$$b_W: \1\to W\otimes W^*\quad {\rm and} \quad d'_W: W\otimes W^*\to \1$$
such that for $a\in \End(W)$ we define
\begin{equation}\label{Trqdef}
Tr_q(a)=d_W'(a\otimes 1)b_W.
\end{equation}
We remark that while $Tr_q$ does not satisfy the trace property $Tr_q(ab)=Tr_q(ba)$ in general,
its restriction to $\End_{U_q\sl_N}(W)$ is indeed a trace. There exists an element $q^{2\rho}\in U_q\sl_N$
such that $Tr_q(a)=Tr(aq^{2\rho})$. We will only need to know its action on $W=V^{\otimes n}$, where we have
\begin{equation}\label{Markovdensity}
Tr_q(a)\ =\ Tr(aD^{\otimes n}),\quad a\in \End(V^{\otimes n}),\quad D = diag(q^{2i-N-1}).
\end{equation}
More generally, we can define a partial trace (also referred to as a contraction, or a conditional expectation)
$\E_X: \End(X\otimes W)\to \End(X)$  by
\begin{equation}\label{condexp}
\E_X(a)\ =\ (1_X\otimes d_W')(a\otimes 1_{W^*})(1_X\otimes b_W),
\end{equation}
which satisfies $Tr_q(a)=Tr_q(\E_X(a))$.
If we take $X=W=V=\C^N$, the vector representation of $U_q\sl_N$, it follows from the discussion above
that $\E_X(U)=Tr_q(U)1_V$, as $V$ is irreducible. One deduces from this more generally that
\begin{equation}\label{Markovcat}
\E_n(U_n)\ =\ Tr(U_n)1,
\end{equation}
where $\E_n$ is the partial trace from $\End(V^{\otimes n+1})$ to $\End(V^{\otimes n})$.

As usual, we define the $q$-number $[k]=(q^k-q^{-k})/(q-q^{-1})$. Then the $q$-dimension of the irreducible 
$U_q\sl_N$  module $V_\mu$ with highest weight $\mu$ is given by
\begin{equation}\label{U(N)dimension} 
\dim_q V_\mu\ =\ \prod_{1\leq i<j\leq N} \frac{[\mu_i-\mu_j+j-i]}{[j-i]}.
\end{equation}
The dimension of a simple $U_q\sp_{2k}$-module $V_\la$ labeled by the Young diagram $\la$ is given by
\begin{equation}\label{Spdimension}
\dim_q V_\la\ =\ \prod_{i<j}\frac{[\la_i-\la_j+j-i][ \la_i+\la_j+2k+2-i-j]}{[j-i][2k+2-i-j]}\ \prod_{i=1}^k \frac{[2\la_i+2k+2-2i]}{[2k+2-2i]}.
\end{equation}

\begin{remark}\label{MarkovHecke} Let $\HnN$ be the image of the Hecke algebra $H_n$ in the representation in $V^{\otimes n}$.
Then we refer to the trace $tr$ on $H_n$ such that $tr(p_\mu)=\dim_qV_\mu/[N]^n$ (as in \ref{U(N)dimension})
 for a minimal idempotent $p_\mu$ in the direct summand of
$H_n$ labeled by $\mu$ as the Markov trace on $\HnN$. It was shown in \cite{WHecke} that the quotient
$\overline{H}_n$ modulo the annihilator ideal of $tr$ is indeed isomorphic to $\HnN$ for $q$ not a root of unity.
\end{remark}

\subsection{Posing the question} We can now make the above mentioned problem of finding a $q$-deformation
for the embedding $Sp(N-1)\subset SL(N)$ more precise as follows: 
We fix $N$  and denote by $\HnN$ the image of the Hecke algebra $H_n$ in the representation in $V^{\otimes n}$.
Then we want to find extensions $C_n=C_N^{(N)}$ of $\HnN$ such that

(a) $C_n\subset C_{n+1}$,

(b) $\HnN\subset C_n\cong  \End_{SL(N)}(V^{\otimes})\subset \End_{Sp(N-1)}(V^{\otimes})$ at least for $q$ not a root of unity,

(c) (Markov property) There exists a normalized trace $tr$ on $C_n$ extending the Markov trace on the Hecke algebra
$\HnN$ compatible with the embedding $C_n\subset C_{n+1}$ for all $n$ such that $tr(e)=1/[N]$ and
\begin{equation}\label{Markovprop}
tr(cg_n)=tr(c)tr(g_n)\quad{\rm for\ all\ }c\in C_n.
\end{equation}
The following lemma will only be needed for $n\leq 2$ in this paper (but more later).

\begin{lemma}\label{tracerestriction} If a trace $tr$ on $C_n$ as in (c) exists, its value on a minimal idempotent $p_\la$
in the direct summand of $C_n$ labeled by $\la$ is equal 
to $\dim_q V_\la/[N]^n$, with $\dim_qV_\la$ given in $\ref{Spdimension}$.
\end{lemma}

$Proof$. We view $D$ as the representation of an element of $Sp(N-1)$ with eigenvalues $q^{\pm 2j}$, $1\leq j\leq (N-1)/2$. 
These are exactly the eigenvalues of the element $q^{2\rho'}$ for $\rho'$ half the sum of the roots of $Sp(N-1)$
in its vector representation $V'$. Hence the value of $tr$ for $\tilde p_\la$ is given
by the symplectic character $\chi^\la(q^{2\rho'})/[N]^n$, where $\tilde p_\la$ is a minimal idempotent in $\End_{Sp(N-1)}(V^{\otimes n})$.
It is well-known that an irreducible $SL(N)$-module labeled by $\mu$  decomposes as a direct sum of $Sp(N-1)$-modules 
labeled by Young diagrams with fewer boxes than $\mu$ except possibly for $\mu$ itself (see  \cite{Brauer}). By assumption $H_n\subset C_n$
has the same inclusion pattern as $\End_{SL(N)}(V^{\otimes n})\subset \End_{Sp(N-1)}(V^{\otimes n})$.
Hence $tr(p_\la)$ is already determined by the values $tr(p_\nu)$ for minimal
idempotents $p_\nu\in C_n$ labeled by diagrams $\nu$ with $|\nu|<|\la|$ and by
the value of $tr$ on minimal idempotents of $H_n$. Hence we can show the claim by induction on $|\la|$.

\subsection{First structure coefficients}
As a simple but important special case of the tensor product rule \ref{tensor}, we obtain that
 $V^{\otimes 2}$ decomposes as an $Sp(N-1)$-module into the direct sum
\begin{equation}\label{secondtensor}
V^{\otimes 2}\ \cong\ 2 \cdot \1 \oplus 2V'\oplus V_{2\Lambda_1}\oplus V_{\Lambda_2};
\end{equation}
here $\1$ is the trivial representation and $2\Lambda_1$ and $\Lambda_2$ denote the Young diagrams with two boxes
in the same row respectively in the same column; see also the Bratteli diagram in Section \ref{Bratteli}. The labeling is  chosen such that $V_{\Lambda_2}$ is 
the nontrivial  irreducible representation appearing in the antisymmetrization of $V^{\otimes 2}$.
It follows from the tensor product rules \ref{tensor} that exactly one copy of $\1$ and of $V$ appears in the symmetrization of $V$, and hence also exactly one
copy of each of these representations appears in the antisymmetrization of $V^{\otimes 2}$.
 It follows from \ref{secondtensor} that $C_2$ has two 
nonisomorphic 2-dimensional irreducible representations, and two one-dimensional representations.
If $e$ denotes the projection onto the trivial representation,
we can choose bases for the two 2-dimensonal representations for which the elements $e$ and $u$ are given by
\begin{equation}\label{erep}
e\ \mapsto\ \left[\begin{matrix} 1&0\\0&0 \end{matrix}\right]\ \oplus\ \left[\begin{matrix} 1&0\\0&0 \end{matrix}\right],
\end{equation}
\begin{equation}\label{urep}
u_1\ \mapsto\ \left[\begin{matrix} a&\sqrt{ab}\\\sqrt{ab}&b \end{matrix}\right]\ \oplus\ \left[\begin{matrix} c&\sqrt{cd}\\\sqrt{cd}&d \end{matrix}\right].
\end{equation}
Using \ref{tensor} again, we obtain (for $N\geq 5$) that
\begin{equation}\label{thirdtensor}
V^{\otimes 3}\ \cong\ 4\cdot \1 \oplus 6V\oplus ... .
\end{equation}
For $N=3$, we have to make the following adjustments: In Eq \ref{secondtensor}, the representation $V_{\La_2}$ does not appear, 
and in Eq \ref{thirdtensor} the representation $V$ only appears with multiplicity 5.

\begin{lemma}\label{acsquare} The matrix entries above satisfy  $a+b=[2]=c+d$ and, if $ab\neq 0$, $(a-c)^2=1$. Moreover, 
there exists a 4-dimensional representation of the Hecke algebra $H_3$ with $u_1$ given by the matrices in \ref{urep} 
and $u_2$ by the same matrix blocks, with the second and third basis
vectors permuted.
\end{lemma}

$Proof.$ The first statement follows from the fact that each of the matrix blocks in \ref{urep} has rank 1 and its only possible nonzero eigenvalue is $[2]$.
We now consider the representation of the Hecke algebra $H_3$ on the four-dimensional space $\Hom(\1,V^{\otimes 3})$. 
It follows from the definition of path bases that $u_1$ can be represented by the 4 by 4 matrix $\rho(u_1)$ with the two diagonal blocks as in Eq \ref{urep}. As $\C v_0\otimes V^{\otimes 2}\cong V^{\otimes 2}$ as an $Sp(N-1)$ module, we can also
assume that $u_2$ is represented by the first matrix block in \ref{urep} on the span of the first and third path,
while it is given by a 2 by 2 block with diagonal entries $c'$ and $d'$ on the remaining two paths.
Checking the Hecke algebra relation \ref{Hecke3} for the (1,1) entry, we deduce that $c'=c$, and hence also $d'=d$.
Checking relation \ref{Hecke3} for entries (1,2) and (2,1), one deduces that the off-diagonal entries in the $c'$,$d'$ block of $u_2$ are equal to $\sqrt{cd}$. This determines the matrix for $u_2$ as claimed.
But then the (21) entry of both the left and the right hand side of  relation \ref{Hecke3} reads as
$$\sqrt{ab}(a^2+bc-1)\ =\ \sqrt{ab}(ac+cd).$$
Substituting $b=[2]-a$ and $d=[2]-c$, and dividing by $\sqrt{ab}$, we obtain
$(a-c)^2=1$, as claimed. 

\begin{corollary}\label{addrelation}
Let $\eth$ be the projection onto the path $0\to 0\to 0\to 0$. 

(a) If $a=c+1$, then $(u_1u_2-u_1)\eth = (u_2u_1-u_2)\eth$, or, equivalently, $(g_1g_2+g_1)e_{(3)}=(g_2g_1+g_2)e_{(3)}$,

(b)  If $a=c+1$, then $(u_1u_2+u_1)\eth = (u_2u_1+u_2)\eth$, or, equivalently, $(g_1g_2-g_1)e_{(3)}=(g_2g_1-g_2)e_{(3)}$.
\end{corollary}

$Proof.$ This follows from the explicit matrices as described in the proof of Lemma \ref{acsquare}.

\subsection{Relations from Markov property} Recall that the Markov property \ref{Markovprop}
 requires $tr(cu_n)=tr(c)tr(u_n)$ for any $c\in C_n$.

\begin{lemma}\label{relation1} 
The Markov property only holds if the matrix entries in \ref{urep} are as follows (where $N=2k+1$):

$$a\ =\ \frac{[k](q^{1/2}+q^{-1/2})}{[k+1/2]},\hskip 6em c\ =\ \frac{[k-1/2]}{[k+1/2]},\hskip 5.5em {\rm if\ } a=c+1,$$

$$a\ =\ \frac{-(q^{1/2}-q^{-1/2})(q^k-q^{-k})}{q^{k+1/2}+q^{-k-1/2}},\hskip 2em c\ =\ \frac{q^{k-1/2}+q^{-k+1/2}}{q^{k+1/2}+q^{-k-1/2}},\hskip 2em {\rm if\ } a=c-1.$$
\end{lemma}

$Proof.$ 
It follows from the definitions and Lemma \ref{tracerestricction} that
$tr(e)=1/[N]$ and 
$$tr(u)\ =[2]\ \frac{[N][N-1]}{[2]}\frac{1}{[N]^2}\ =\ \frac{[N-1]}{[N]}.$$
Using the explicit matrix representations \ref{erep} and \ref{urep} and the weights of the traces in Lemma \ref{tracerestriction}, we obtain
$$tr(ue)\ =\ \frac{1}{[N]^2}(a + c([N]-1)).$$
It follows from the Markov property $tr(ue)=tr(u)tr(e)$, with $N=2k+1$ that
$$a+c([2k+1]-1)=[2k].$$
By Lemma \ref{acsquare}, we have $a=c\pm 1$. Let us consider the case $a=c+1$.
Substituting this into the last equation we obtain
$$c=\frac{[2k]-1}{[2k+1]}\ =\ \frac{[k-1/2]}{[k+1/2]}\ =\ \frac{q^{k-1/2}-q^{-k+1/2}}{q^{k+1/2}-q^{-k-1/2}}.$$
The formula for $a$ follows from $a=c+1$. The case $a=c-1$ goes similarly.

\begin{remark} \label{modrestriction} The previous lemma implies that there are at most two module categories $\M$ of $Rep(U_q\sl_N)$ for our given
fusion rules if they allow categorical traces which also satisfy the compatibility condition $Tr_M(\E_M(a))=Tr_{M\otimes W}(a)$ for
$a\in \End_\M(M\otimes W)$, where $\E_M$ is defined for an object $M$ in $\M$ and an object $W\in Rep(U_q\sl_N)$ as in \ref{condexp}.
\end{remark}

\subsection{Correction for $q\to 1$} We assume in this section that  $a=c-1$. It follows from Lemma \ref{relation1} that  $ab=0$ for $q=1$ in this case. This would make the
representation in the first matrix block in \ref{urep} reducible. This can be avoided by introducing the elements
$u_{12}$ and $u_{21}$ defined below; they correspond to intertwiners between the two copies of the trivial representation
in $V^{\otimes 2}$.  We shall show later that one can make sense of them also at $q=1$.
 It will be convenient to define
\begin{equation}\label{Npdef}
\Np=\frac{q^{N/2}-q^{-N/2}}{q^{1/2}-q^{-1/2}}\quad {\rm and} \quad
\Nm=\frac{q^{N/2}+q^{-N/2}}{q^{1/2}+q^{-1/2}}.
\end{equation}
Observe that
$$b\ =\ [2]-a\ =\ \frac{q^{k+1}+q^{-k-1}}{[N]_-},$$
$$ \sqrt{ab}=\frac{-i(q^{1/2}-q^{-1/2})}{[N]_-}\sqrt{[N]-1},$$
where we used the identity $[k](q^{k+1}+q^{-k-1})=[N]-1$, and where the choice of sign for $\sqrt{ab}$
will turn out to be immaterial. 
Let $e_{(2)}$ be the subprojection of $e$ which is nonzero only in the first matrix block in \ref{erep}. We define the element
\begin{align}\label{u21def}
u_{21}\  &=\  \frac{[N]_-}{-i(q^{1/2}-q^{-1/2})}(1-e_{(2)})u_1e_{(2)}\\
 &=\ \frac{i[N]_-}{-q^{1/2}-q^{-1/2}}u_1e_{(2)}+i(q^{1/2}-q^{-1/2})[k] e_{(2)},\notag
\end{align}
where we used $e_{(2)}u_1e_{(2)}=ae_{(2)}$. We similarly define $u_{12}=u_{21}^t$ by the same expression as in \ref{u21def}
with $u_1$ and $e_{(2)}$ interchanged. By construction, it follows that the elements $u_{21}$ and $u_{12}$ are represented
by the matrices
$u_{21} =\sqrt{[N]-1}E_{21}$ and $u_{12}=\sqrt{[N]-1}E_{12}$, where $E_{ij}$ are matrix units in the first $2\times 2$ matrix block in \ref{urep}.
The following results follow immediately from this.

\begin{lemma}\label{u21lemma}
We have $u_{12}u_{21}=([N]-1)e_{(2)}$, $u_{21}u_{12}e_{2)}=0=e_{(2)}u_{21}u_{12}$ and the element $P=e_{(2)}+u_{12} + u_{21}+ u_{21}u_{12}$
satisfies $P^2=[N]P$.
\end{lemma}

\subsection{Summary of relations} The following is a preliminary definition of the algebras $C_n$. 
The precise, but less intuitive definition of the algebras $C_n$ in various versions will be given
in the following section.

\begin{definition}\label{allrelations} Fix $N=2k+1$. Then we define the algebra $C_{n,\pm}=C^{(N)}_{n,\pm}$ via generators
$e$, $u_i$, $1\leq i<n$ with the following relations:

(a) The elements $u_i$ satisfy the Hecke algebra relations,

(b) We have a sequence of idempotents $e_{(r)}$ defined inductively by  $e_{(0)}=1$, $e_{(1)}=e$ and
$$e_{(r+1)}\ =\ e_{(r)}u_re_{(r)} - \frac{q^{(N-2)/2}-q^{-(N-2)/2}}{q^{N/2}-q^{-N/2}}e_{(r)},\quad {\rm for\ } C_{n,+},$$
$$e_{(r+1)}\ =\ \frac{q^{(N-2)/2}+q^{-(N-2)/2}}{q^{N/2}+q^{-N/2}}e_{(r)}-e_{(r)}u_re_{(r)},\quad {\rm for\ } C_{n,-}.$$

(c)   For $j<r$ we have $(u_{j-1}u_j-u_{j-1})e_{(r)} = (u_{j}u_{j-1}-u_j) e_{(r)}$ for $C_{n,+}$ and  $(u_{j-1}u_j+u_{j-1}) e_{(r)} = (u_2u_1+u_2) e_{(r)}$ for $C_{n,-}$.
\end{definition}

\begin{remark}\label{relationsremark}
1. In spite of fractions in the exponents, it is not hard to check that the relations only depend on $q$ and not on the choice
of a square root $q^{1/2}$. E.g. we have
$$\frac{q^{(N-2)/2}-q^{-(N-2)/2}}{q^{N/2}-q^{-N/2}}\ =\ \frac{q^{N-1}+q^{N-2}+\ ...\ +q^{1-N}}{q^N+q^{N-1}+\ ...\ +q^{-N}}.$$

2. Let $e_{\pm}(q)$ and $u_{i,\pm}(q)$ be the generators  of $C_{n,\pm}(q)$ for a given choice of $q$ respectively. 
Then we can check that the maps
$$e_-(q)\mapsto e_+(-q), \hskip 3em u_{i,-}(q)\mapsto -u_{i,+}(-q),\quad 1\leq i<n,$$
define an isomorphism between $C_{n,-}(q)$ and $C_{n,+}(-q)$.
\end{remark}

\subsection{Alternative definitions of the algebras $C_n$} We make the following adjustments for the precise definition
of the algebra $C_n=C_n(p,q)$ as an algebra depending on two variables $p$ and $q$. First, we substitute $p=q^N$ in the previous relations.
Secondly, we introduce additional generators which will only be relevant for the important classical limit $q\to 1$
for $C_{n,-}$ (respectively for $q\to -1$ for $C_{n,+}$).

\begin{definition}\label{pqdef} (Two variable definition) We define the algebra $C_n=C_{n,\pm}=C_{n,\pm}(p,q)$ 
over the field $\C(p,q)$ of rational functions
in the variables $p$ and $q$ as follows: We have generators $u_i$, $e_{(r)}$ for $1\leq i,r<n$
with the following relations:

(a) The elements $u_i$ satisfy the Hecke algebra relations,

(b) The  elements $e_{(r)}$ are idempotents which satisfy the relations
$$e_{(r+1)}\ =\ e_{(r)}u_re_{(r)} - \frac{p^{1/2}q^{-1}-qp^{-1/2}}{p^{1/2}-p^{-1/2}}e_{(r)},\quad {\rm for\ } C_{n,+},$$
$$e_{(r+1)}\ =\ \frac{p^{1/2}q^{-1}+qp^{-1/2}}{p^{1/2}+p^{-1/2}}e_{(r)}-e_{(r)}u_re_{(r)},\quad {\rm for\ } C_{n,-}.$$

(c)   For $j<r$ we have $(u_{j-1}u_j-u_{j-1})e_{(r)} = (u_{j}u_{j-1}-u_j) e_{(r)}$ for $C_{n,+}$ and  $(u_{j-1}u_j+u_{j-1}) e_{(r)} = (u_2u_1+u_2) e_{(r)}$ for $C_{n,-}$.
\end{definition}
\begin{remark} Similarly as e.g. for the algebras defined in \cite{BBMW}, this definition is not convenient if we are interested in
obtaining the classical limits for E$q\to 1$ and $p=q^N\to 1$. This can be addressed by introducing additional generators.
It will be shown in Section \ref{tensorsection} that one can make sense of these additional elements if $q\to 1$.
\end{remark}

\begin{definition}\label{Zdef} (Extended definition for $C_n(q^N,q)$) We now let $p=q^N$ as before. 
We add to the usual generators $u_i$, $ e_{(r)}$ for $1\leq i,r<n$ also the elements
$$\frac{1}{q^{1/2}+q^{-1/2}}u_i e_{(r)},\ \frac{1}{q^{1/2}+q^{-1/2}} e_{(r)}u_i,\quad 1\leq i<r<n,\quad  {\rm for \ }C_{n,+},$$
$$\frac{1}{q^{1/2}-q^{-1/2}}u_i e_{(r)},\ \frac{1}{q^{1/2}-q^{-1/2}} e_{(r)}u_i,\quad 1\leq i<r<n,\quad  {\rm for \ }C_{n,-}.$$
The relations are the same as in Definition \ref{allrelations}.
\end{definition}

\subsection{Basic structure results} We will prove existence of nontrivial representations of the algebras $C_{n,\pm}$ in the next section.

\begin{proposition}\label{basicstructure} (a) The map $\Phi$ given by $e\mapsto e_{(r+1)}$, $u_i\mapsto e_{(r)}u_{r+i}$ defines a homomorphism of $C_{n-r}$ onto $e_{(r)}C_ne_{(r)}$.

(b) The algebras $C_{n,\pm}$ are spanned by the $H_n-H_n$ bimodules $H_ne_{(r)}H_n$, $0\leq r\leq n$. In particular, they are finite dimensional.

(c) The span $I_r$ of $\cup_{s\geq r} H_ne_{(s)}H_n$ is a two-sided ideal of $C_n$ for $1\leq s\leq n$.
\end{proposition}

$Proof.$ The homomorphism property in (a) follows directly from the relations. For (b), we observe that the claimed spanning set 
contains the generators of $C_n$. It hence suffices to show that multiplying it by a generator from the right or left
will still produce an element in the span. This is obviously true for the generators $u_i$.
We prove the claim for multiplication by $e$  by induction on $n$, with the
statement obviously true for $n=1$. Observe that the Hecke algebra $H_n$ is spanned by elements of the form
$au_1b$ or $a$, with $a,b\in H_{2,n}$, where $H_{2,n}$ is the subalgebra generated by $u_2, u_3,\ ...\ u_{n-1}$ and 1.
We then have, using $e_{(r)}=ee_{(r)}$
$$e(au_1b)e_{(r)}=a(eu_1e)be_{(r)}= c\ abe_{(r)} + ae_{(2)}be_{(r)}.$$
But now $e_{(2)}be_{(r)}\in e_{(2)}H_{2,n}e_{(2)}\subset \Phi(eH_{n-1}e)$ by (a). Hence, by induction assumption, we have
$$e_{(2)}be_{(r)}\in span(\cup_s \Phi^{-1} (H_{n-1}e_{(s)}H_{n-1})\subset span \cup_s  H_{n-1}e_{(s)}H_{n-1}.$$
The proof for multiplication by $e$ from the right goes completely analogously. This finishes the proof of statement (b).
The surjectivity statement in (a) now follows from (b), as $\Phi$ maps $e_{(s)}$ to $e_{(s+1)}$.
Our proof of (b) also implies statement (c).

\subsection{Dimension estimates}\label{dimestimsec}  We will explicitly construct spanning sets for the algebras $C_n$ which will later  be shown to be bases.
To do so, we shall use two well-known facts about Hecke algebras, here only formulated for
Hecke algebras for type $A$ (see \cite{humphreys} for details). Let $s_i$, $1\leq i<n$ be a set
of simple reflections for the symmetric group $S_n$ (say $s_i=(i,i+1)$). Then any element $w\in S_n$ can be written as a product
of simple reflections. Any such expression for $w$ with the minimal number of factors is called a reduced word,
and the number of factors is called the length $\ell(w)$ of $w$. In the case of the symmetric
group, the length $\ell(w)$ can also be defined as the number of pairs $i<j$ such that $w(i)>w(j)$.
 Replacing the elements $s_i$ by generators $g_i$
in such an expression defines an element $h_w\in H_n$ which does not depend on the choice of reduced expression for $w$.

It is easy to see that the shortest elements in the left cosets of $S_r\subset S_n$ are given by  permutations $w$ which satisfy
$w(i)<w(j)$ for any $1\leq i<j\leq r$, and that the shortest elements in the left cosets of $S_r\times S_{n-r}\subset S_n$
are given by permutations $w$ which
also satisfy the additional conditions $w(i)<w(j)$ for all $r<i<j\leq n$. It is well-known (and easy to check for these
cases)  that each such coset contains a unique element of
lowest length.

\begin{definition}\label{spandef}
If $w\in S_n$ and $h_w$ the corresponding element in $H_n$, we define $h_w^T=h_{w^{-1}}$. If $\S\subset H_n$
we define $\S^T=\{ h^T, h\in \S\}$. Using these conventions, we define 
$\Ca_{n,r}=\{h_w\}$, where $w$ ranges over the shortest elements in the left cosets of $S_r\subset S_n$,
and we define $\D_{n,r}=\{ h_w\}$, where now $w$ ranges over the shortest elements of the left
cosets of  $S_r\times S_{n-r}\subset S_n$. Finally, we define $\B_r$ inductively by  $\B_{0}=\emptyset$, $\B_{1}=\{ 1\}$ and
$$\B_{r+1}\ =\ \B_{r}\  \cup\ \bigcup_{j=1}^{r} g_jg_{j+1}\ ...\ g_{r}\B_{r-1}.$$
\end{definition}
We remark that $|\Ca_{n,r}|=[S_n:S_r]=n!/r!$ and that $|\D_{n,r}|=\binom{n}{r}$, and that $\D_{n,r}^T$
contains the elements $h_w$ with $w$ running through the shortest elements of the right cosets of $S_r\times S_{n-r}\subset S_n$.

\begin{lemma}\label{Bnestimates}
(a) We have $|\B_{r}|=h_r$, with $h_r$ as in \ref{hr}.

(b)  The set $\Ca_{n,r}\B_{r}e_{(r)}$ spans $H_ne_{(r)}$.
\end{lemma}

$Proof.$ Statement (a) follows from the definitions of $h_r$ in \ref{hr} and $\B_{r}$. To establish statement (b),
let us first prove it in the special case $n=r$ by induction on $r$. This is obviously true for $r=1$.
We now prove by downward induction from $s=r-1$ to $s=1$ that
$$H_rg_rg_{r-1}\ ...\ g_se_{(r+1)}\ \subset H_re_{(r+1)} + H_rg_re_{(r+1)}.\eqno(*)$$
For $s=r-1$  it follows from relation (c) that $g_rg_{r-1}e_{(r+1)}=(g_{r-1}g_r-g_{r-1}+g_r)e_{(r+1)}$.
This implies $H_rg_rg_{r-1}e_{(r+1)}$ is contained in the right hand side of the claim.
The induction step for $s<r-1$ is shown in the same way. 
As $H_{r+1}= H_r+\sum_{s=1}^rH_rg_rg_{r-1}\ ...\ g_s$, it follows that 
\begin{align}
H_{r+1}e_{(r+1)}\ &=\ H_re_{(r+1)} \ +\ H_rg_re_{(r+1)}\\ \notag
&=\ H_re_{(r+1)} \ +\ \sum_{j=1}^rH_{r-1}g_{r-1}g_{r-2}\ ... g_je_{(r-1)}g_re_{(r+1)},\notag
\end{align}
where the summand for $j=r$ is defined to be equal to $H_{r-1}e_{(r-1)}g_re_{(r+1)}$. 
The claim for $n=r$ follows from $(*)$ and the definition of $\B_{r+1}$. 
The general claim for $n>r$ follows from this and the fact that $H_n=\Ca_{n,r}H_r = \bigcup_{h_w\in\Ca_{n,r}} h_wH_r$.

\begin{lemma}\label{estimates} We have the following inequalities:

(a) $\dim H_ne_{(r)}\leq \frac{n!}{r!}h_r$,

(b) $\dim H_ne_{(r)}H_n\leq \frac{n!}{r!}\binom{n}{r}h_r^2$,

(c) We have $\dim C_n \leq \sum_{r=0}^n \frac{n!}{r!}\binom{n}{r}h_r^2$.
\end{lemma}

$Proof.$ It follows from the inductive definition of $\Ca_{r,n}$ that $|\Ca_{r,n}|=n!/r!$ 
This and Lemma \ref{Bnestimates} imply claim (a). 
The  inequality in (b) follows from (a)
and the surjective map $H_n\otimes_{H_{r+1,n}}H_n\to H_ne_{(r)}H_n$.  Finally, claim (c) follows from the fact that 
$C_n$ is spanned by the subspaces $H_ne_{(r)}H_n$, $0\leq r\leq n$.

\subsection{An explicit spanning set} We use the notations from Definition \ref{spandef}.

\begin{proposition}\label{spanningset}
The set $\B(r)=\Ca_{r,n}\B_{r,n}e_{(r)}\B_{n,r}^T\D_{r,n}^T$ spans $H_ne_{(r)}H_n$ for $1\leq r\leq n$,
and hence $\B=\bigcup \B(r)$ spans $C_n$.
\end{proposition}

$Proof.$ We hace already proved in Lemma \ref{Bnestimates} that  set $\Ca_{n,r}\B_{r}e_{(r)}$ spans $H_ne_{(r)}$.
One can show the same way that the set $e_{(r)}\B^T(r)$ spans $e_{(r)}H_r$. The claim now follows from this and
the fact that $\D_{n,r}^T$ contains the elments of minimal lengths for all right cosets of $H_r\times H_{n-r}\subset H_n$.

\section{Tensor product representations}\label{tensorsection}

We will give explicit representations of the algebras $C_n=C_{n,+}$ in this section. In view of Remark \ref{relationsremark}, 
these representations can be easily modified to representations of the algebras $C_{n,-}$.  As before, let 
$V=\C^N$ with $N=2k+1$ and let $\{ v_i\}$ denote the standard basis of $\C^N$. We define $v_o=\sum_{i=1}^N \al_iv_i$,
with $\al_i=q^{(k+1-i)/2}/\|v_o\|$, where $\|v_o\|^2=\sum \al_i^2$.

\subsection{A matrix for $e$} 
We define the $N$ by $N$ matrix $E=(e_{ij})$ by $e_{ij}=\al_i\al_j$. Moreover, we modify
the Hecke algebra representation in  \ref{Heckerep}  to $u\ \mapsto \ U$, where the matrix $U$ is defined by
\begin{equation}\label{Heckerep2}
U_{|span\{ v_i\otimes v_j, v_j\otimes v_i\}}\ =\ 
\left[\begin{matrix} q^{-1}&1\\1&q \end{matrix}\right], \quad i<j.
\end{equation}
with the matrices $U_i$  defined as in  \ref{Hecketensor}.
We then consider the map $\Phi$ which maps the generators of $C_n$ into $\End(V^{\otimes n})$ given by
\begin{equation}\label{Phidef}
\Phi: \quad e\mapsto E\otimes 1_{n-1},\hskip 3em u_i\mapsto U_i,\quad 1\leq i<n.
\end{equation}

\begin{lemma}\label{Erep} The map $\Phi$ is compatible with relation (b), mapping the
element $e_{(n)}$ to $E^{\otimes n}$. In particular,
we have 
$$E\otimes E\ =\ (E\otimes 1)U(E\otimes 1)\ -\ \frac{q^{k-1/2}-q^{-k+1/2}}{q^{k+1/2}-q^{-k-1/2}}\ E\otimes 1.$$
\end{lemma}

$Proof.$ 
The vector $v_o$ spans the image of $E$ by definition. We then calculate
$$U(v\otimes v_j)\ =\ \sum_{i=1}^{j-1} \al_i(q^{-1}v_i\otimes v_j +v_j\otimes v_i) \ 
+\ \sum_{i=j+1}^{N} \al_i(q v_i\otimes v_j + v_j\otimes v_i),$$
$$(E\otimes 1)U (v\otimes v_j)\ =\ 
\sum_{i=1}^{j-1} \al_i^2q^{-1}v_i\otimes v_j + \al_i\al_j v_j\otimes v_i \ 
+\ \sum_{i=j+1}^{N} \al_i^2q v_i\otimes v_j + \al_i\al_j v_j\otimes v_i \ =$$
$$=\ \al_j\sum_{i=1}^N v\otimes \al_iv_i \ +\ \beta v\otimes e_j,$$
where
$$\beta\ =\ 
-\al_j^2\ +\ \sum_{i=1}^{j-1} \al_i^2q^{-1}\ 
+\ \sum_{i=j+1}^{N} \al_i^2q.$$
It is now straightforward to check that $\al_j\sum_{i=1}^N v\otimes \al_iv_i = (E\otimes E)(v\otimes v_j)$
and
$$\beta = \frac{q^{k-1}+q^{k-2}+\ ...\ + q^{1-k}}{\| v_o\|^2}\ =\ \frac{q^{k-1/2}-q^{-k+1/2}}{q^{k+1/2}-q^{-k-1/2}},$$
where the last equality is obtained by multiplying both numerator and denominator by $q^{1/2}-q^{-1/2}$.
This proves the second statement in the claim. We can now show by induction on $n$ that $\Phi$ can be extended to a homomorphism which maps $e_{(n)}$ to $E^{\otimes n}$.

\subsection{Checking relation (c)} Let $\{ v_i\}$ and $v_o$ be as above. Then we have

\begin{lemma}\label{relationclemma}
$$(U_1U_2-U_1)v^{\otimes 3}\ =\ (U_2U_1-U_2)v^{\otimes 3}.$$
\end{lemma}

$Proof$. This is a straightforward calculation.  It is manageably tedious if one checks it separately
on the span of all possible permutations of $v_i\otimes v_j\otimes v_m$ for any choice of indices $i,j$ and $m$.
 If all indices are mutually distinct,
one obtains two $6\times 6$ matrices for $U_1$ and $U_2$, both with three $2\times 2$ blocks. 
Moreover, the coefficient for each of these vectors in the basis expansion of $v^{\otimes 3}$ is equal
to $\al_i\al_j\al_m$. Hence it suffices to check the claim for these two $6\times 6$ matrices,
applied to the vector $(1,1,1,1,1,1)^T$, which is not very hard. The case where two indices coincide
is done similarly and easier, only involving $3\times 3$ matrices.

\begin{remark} If we define the algebra $\tilde C_3$ like the algebra $C_3$ without relation (c), it can be shown
that $\tilde C_3$ modulo its radical is isomorphic to $C_3$.
\end{remark}

\subsection{Classical limits}\label{climits} It will be more convenient to consider the representations for $C_{n,-}$,
i.e. we basically replace $q$ by $-q$ in the matrix $\Phi(e)$ and in the coefficients of the vector $v_o$, see Remark \ref{relationsremark}.
We are going to show that the elements $\frac{1}{q^{1/2}-q^{-1/2}}u_ie_{(r)}$ still make sense in our representation
even at $q=1$.

\begin{lemma}\label{limitlemma} (a) The matrix coefficients of  $\frac{1}{q^{1/2}-q^{-1/2}}U_iE_{(r)}$ in the tensor product representation
of $C_{n,-}$ are well-defined also at $q=1$, up to the choice of the square root $q^{1/2}$.

(b) The elements $E^{\otimes 2}$, $\frac{1}{q^{1/2}-q^{-1/2}}U_1E^{\otimes 2}$ and
$\frac{1}{q^{1/2}-q^{-1/2}}E^{\otimes 2}U_1$  generate an algebra which is isomorphic to the $2\times 2$ matrices
if $[N]\neq 1$.
\end{lemma}

$Proof.$ We calculate
$$U\vo^{\otimes 2}\ =\ \sum_{i<j}\al_i\al_j(q^{-1}-1)v_i\otimes v_j+(q-1)v_j\otimes v_i).$$
It follows that also the coefficients in the expression for $\frac{1}{q^{1/2}-q^{-1/2}}U\vo^{\otimes 2}$ are in
$\Z[q^{\pm 1/2}]$. As all columns of $E^{\otimes 2}$ are proportional to $\vo^{\otimes 2}$,  claim (a) follows.

Recall the definition of $u_{21}$ in Definition \ref{u21def}.
It follows from part (a) and the equations below \ref{Npdef} that $\Phi(u_{21})$ and $\Phi(u_{12})$ are well-defined and that
$$\Phi(u_{21}u_{12})\ =\ \frac{-\Nm^2}{(q^{1/2}-q^{-1/2})^2} (1-E^{\otimes 2})UE^{\otimes 2}U (1-E^{\otimes 2})$$
has the nonzero eigenvalue
$$\frac{-\Nm^2 a_-([2]-a_-)}{(q^{1/2}-q^{-1/2})^2}\ =\ [k](q^{k+1}+q^{-k-1})\ =\ [N]-1.$$

\subsection{An embedding of $Sp(N-1)$ into $Sl(N)$} Let $A$ be the $N\times N$ matrix (with $N=2k+1$ odd)
defined by 
$$a_{ij}=\begin{cases} (-q^{-1})^{k+1-(i+j)/2} & if\ i<j,\cr - (-q^{-1})^{k+1-(i+j)/2}  & if\ i>j.\end{cases}$$
As usual, we assume  fixed choices of $(-q)^{1/2}$  and of $(-q^{-1})^{1/2}$ in all these formulas such that their product
is equal to $-1$.

\begin{lemma}\label{skewlemma}
The matrix $A$ has rank $N-1$ in a neighborhood of $q=1$, with kernel $\vo = \sum_{j=1}^{2k+1} (-q)^{(k+1-j)/2}v_j$.
Hence we obtain a symplectic form $(v,w)=v^TAw$  for $q=1$ whose restriction to any complement $V'$ of $\vo$ is nondegenerate.
\end{lemma}

$Proof.$ We check that
\begin{align}
(A\vo)_i\ &=\ \sum_{j=1}^{i-1} -  (-q^{-1})^{k+1-(i+j)/2} (-q)^{(k+1-j)/2}
\ +\ \sum_{j=i+1}^{2k+1}   (-q^{-1})^{k+1-(i+j)/2} (-q)^{(k+1-j)/2}\notag\\
&= \ \sum_{j=2}^{2k+1} (-q^{-1})^{(k+1-i)/2} (-1)^{k+1-j}\ =\ 0.\notag
\end{align}
To determine the rank at $q=1$, we observe that after conjugation by the diagonal matrix $D=diag((-1)^{(k+1-i)/2})$
the matrix entries become equal to $a_{ij}=(-1)^{k+1-i}$ for $i<j$ and $a_{ij}=-a_{ji}$ for $i>j$.
It is now easy to see that the transformed matrix has eigenvectors $(1,0,0,..., \pm 1)$, $(0,1,0,...,\pm 1,0)$ etc.
This proves the claim about the rank.
 
\begin{proposition}\label{classlimprop} Fix $N=2k+1$ and let $V=\C^N$. If $q=1$, the representation of $C_{n,-}$
into $\End(V^{\otimes n})$ surjects onto $\End_{Sp(N-1)}(V^{\otimes n})$, where the embedding of $Sp(N-1)\subset Sl(N)$
is defined via the symplectic form given by the matrix $A$ in Lemma \ref{skewlemma} at $q=1$. In particular, 
 $\dim_{\C(q)} C_{n,-}(q^N,q)\geq
\dim \End_{Sp(N-1)}(V^{\otimes n})$.
\end{proposition}

$Proof.$ The image of $C_{n,-}$ in $\End(V^{\otimes n})$ at $q=1$ contains the usual action of the symmetric group
$S_n$ on $V^{\otimes n}$ and the projection $E\in\End(V)$ onto $\C \vo$. By Lemma \ref{limitlemma}, it also 
acts as a full $2\times 2$ matrix algebra on $\Hom(\1, V^{\otimes 2})$. So, in particular, it also must contain the projection $F$.
It follows from Theorem \ref{fundamental} that $C_{n,-}(q)$ maps surjectively onto $\End_{Sp(N-1)}(V^{\otimes n})$ for $q=1$.
This implies the estimate about the dimensions.

\subsection{A basis for $C_n(p,q)$} It will be convenient to consider the 2-variable version of $C_n$, as defined in 
Definition \ref{pqdef}.

\begin{theorem}\label{dimensiontheorem} The spanning set in Proposition \ref{spanningset} is a basis for the two-variable
version $C_n(p,q)$, viewed as an algebra  over the field of rational functions in $p$ and $q$. In particular, we have
$\dim C_n(p,q)= \sum_{r=0}^n h_r^2\frac{n!}{r!}\binom{n}{r}$. 
\end{theorem}

$Proof.$
It follows from  Proposition \ref{classlimprop} and Lemma \ref{estimates},(c)
that $\dim C_{n,-}(q^N,q)$ is equal to $\dim \End_{Sp(N-1)}(V^{\otimes n})$  if $N>2n$.
 Hence the spanning set in Proposition \ref{spanningset} is a basis for these values. 
In particular, this is true for the two-variable version of $C_{n,-}(p,q)$ if $p=q^N$ with $N>2n$ odd.
Hence it is true in general by Zariski density. The claim can be similarly shown for $C_{n,+}(p,q)$ using
Remark \ref{relationsremark}.

\section{Markov traces}

Recall that we defined $Tr_q$ on $\End(V^{\otimes n})$ in Section \ref{tracesec} by $Tr_q(a)=Tr(aD^{\otimes n})$,
where $a\in \End(V^{\otimes n})$ and $D=diag(q^{2i-N-1})$, see \ref{Trqdef} and also \ref{Markovdensity}. 
We define a functional $\phi$ on $C_n=\CnN$ as a normalized pull-back of $Tr_q$ by

\begin{equation}\label{Markovtrdef}
\phi(c)\ =\ \frac{1}{[N]^n} Tr_q(\Phi(c)),\quad c\in C_n.
\end{equation}
Observe that $\phi(1)=1$. It is the goal of this section to show that the functional $\phi$ defines a trace $tr$ on $C_n$
which satisfies the Markov  condition \ref{Markovprop}.

\begin{lemma}\label{Markovlemma}
The functional $\phi$ has the following properties.

(a) The restriction of $\phi$ to $H_n$ defines a trace.

(b) $\phi(cg_{n-1})=\phi(c)\phi(g_{n-1})$ for all $c\in C_{n-1}$.

(c) $\phi(e_{(r)})=\phi(e)^r=\frac{1}{[N]^r}$.

(d) $\phi(ch)=\phi(hc)$ for all $h\in H_n$, $c\in C_n$
\end{lemma}

$Proof.$ Part (a) follows from the  discussion in Section \ref{tracesec}, and part (b) follows from \ref{Markovcat}.
We also have for $n=r=1$ and for $C_{n,+}$
$$Tr_q(E)=Tr(ED)\ =\ \sum_{i=1}^N q^{2i-N-1}q^{(N+1)/2-i}/{[N]_+}\ =\ \sum_{i=1}^N q^{i-(N+1)/2}/{[N]_+}\ =\ 1.$$
One deduces from this that 
$$Tr_q(E^{\otimes r})\ =\ Tr((ED)^{\otimes r})Tr(D^{\otimes n-r})\ =\ [N]^{n-r},$$
from which follows  claim (c).
As $\Phi(h)$ commutes with $D^{\otimes n}$ for all $h\in H_n$, we have
$$Tr_q(\Phi(ch))=Tr(\Phi(c)\Phi(h)D^{\otimes n})=Tr_q(\Phi(c)D^{\otimes n}\Phi(h))=Tr_q(\Phi(h)\Phi(c)D^{\otimes n})
=Tr_q(\Phi(hc)).$$

\subsection{Technical lemmas} It follows directly from the relations that the map 
\begin{equation}\label{Thetadef}
\Theta_n: g_i\mapsto g_{n-i},\quad 1\leq i<n
\end{equation}
induces an automorphism of the Hecke algebra $H_n$ which will be denoted by the same letter.
Also observe that if $w_o\in S_n$ is defined by $w_o(i)=n-i$, and we denote the corresponding
map on $V^{\otimes n}$ given via permutation of the factors by the same letter, we have
\begin{equation}\label{woreversal}
w_oU_iw_o=U_{n-i}(q^{-1}),
\end{equation}
where $U_i(q^{-1})$ is given by the same matrix as $U_i$, with every occurrence of  $q$ replaced by $q^{-1}$.

\begin{lemma}\label{Markovlemm1}
We have
$$Tr_q(\Phi(h_1)E^{\otimes n}\Phi(h_2))=Tr(\Phi(h_1)E^{\otimes n}(q^{-1})\Phi(h_2))=Tr(\Phi(\Theta_n(h_1))E^{\otimes n}\Phi(\Theta_n(h_2))).$$
\end{lemma}

$Proof.$ To avoid cumbersome notation, we denote $H_i=\Phi(h_i)$ in this proof. Then we have
$$Tr_q(H_1E^{\otimes n}H_2)\ =\ Tr(H_1(D^{1/2}ED^{1/2})^{\otimes n}H_2)\ =\ 
Tr(H_1E^{\otimes n}(q^{-1})H_2),$$
from which follows the first equality in the statement. Now observe that the structure coefficients in the 
defining relations of $C_n$ are invariant under $q\leftrightarrow q^{-1}$.
Hence we also obtain a representation of $C_n$ via the assignment
$$e\mapsto E(q^{-1}),\hskip 3em u_i\mapsto U_i(q^{-1}),\quad 1\leq i<n.$$
As $E^{\otimes n}HE^{\otimes n}=\al E^{\otimes n}$ for some scalar $\al$, it also follows 
$E^{\otimes n}(q^{-1})H(q^{-1})E^{\otimes n}(q^{-1})=\al E^{\otimes n}(q^{-1})$ for the same scalar.
Applying this to $H=H_2H_1$, we obtain
$$Tr(H_1E^{\otimes n}H_2)= Tr(E^{\otimes}H_2H_1E^{\otimes n}) =\al = $$
$$= Tr(E^{\otimes}(q^{-1})H_2H_1(q^{-1})E^{\otimes n}(q^{-1}))\ =\ Tr(\Theta_n(H_1)E^{\otimes}(q^{-1})\Theta_n(H_2)).$$
Replacing $H_i$ by $\Theta(H_i)$ for $i=1,2$ in the equation above now proves the second equality in the statement.

\begin{lemma}\label{Markovlemma2} We have $e_{(n)}he_{(n)}=e_{(n)}\Theta_n(h)e_{(n)}$
for all $h\in H_n$.
\end{lemma}

$Proof.$ The claim is proved by induction on $n$, with $n=1,2$ being trivially true.
For the induction step from $n-1$ to $n$ first observe that for any $A\in\End(V^{\otimes n-1})$ we have
$$E^{\otimes n}(A\otimes 1)E^{\otimes n} = E^{\otimes n-1}AE^{\otimes n-1}\otimes E
= E\otimes E^{\otimes n-1}AE^{\otimes n-1}= E^{\otimes n}(1\otimes A)E^{\otimes n},$$
as $ E^{\otimes n-1}AE^{\otimes n-1}$ is a scalar multiple of $E^{\otimes n-1}$. 
We define the homomorphism $sh:H_{n-1}\to H_n$ via $sh(u_i)=u_{i+1}$. If $H\in \Phi(H_{n-1})$, it follows
that $E^{\otimes n} sh(H) E^{\otimes n} =E^{\otimes n} HE^{\otimes n} $. Moreover, by induction assumption,
we have 
$$E^{\otimes n} (\Theta_{n-1}(H)\otimes 1) E^{\otimes n} =E^{\otimes n} (H\otimes 1) E^{\otimes n} .$$ 
Hence we have
\begin{equation}\label{Hncase}
e_{(n)}he_{(n)} =e_{(n)}sh(\Theta_{n-1}(h))e_{(n)}=e_{(n)}\Theta_n(h)e_{(n)},
\end{equation}
 which proves the claim for $h\in H_{n-1}$. Let now $h=h'g_{n-1}$ with $h'\in H_{n-1}$. 
We first observe that for $h\in H_n$, $e_{(n-1)}he_{(n-1)}$ is a linear combination
of $\al e_{(n-1)}+\beta e_{(n)}$. Using $e_{(n)}e_{(n-1)}=e_{(n)}$, one deduces easily that
$$e_{(n)}he_{(n-1)}=e_{(n)}he_{(n)}.$$
If $e_{(n-1)}h'e_{(n-1)}=\gamma e_{(n-1)}$, we calculate
$$e_{(n)}he_{(n)}= e_{(n)}h'e_{(n-2)}g_{n-1}e_{(n)}=\gamma e_{(n)}g_{n-1}e_{(n)}=\gamma (c'+1)e_{(n)}.$$
On the other hand,
$$e_{(n)}\Theta_n(h)e_{(n)} = e_{(n)}\Theta_n(h')g_1e_{(n)}= e_{(n)}\Theta_n(h')eg_1e_{(n)},$$
where we used that $\Theta_n(h')\in H_{2,n}$ commutes with $e$.
But then
$$e_{(n)}\Theta_n(h)e_{(n)} = e_{(n)}\Theta_n(h')(c'+1)e_{(n)}=\gamma (c'+1)e_{(n)},$$
by \ref{Hncase}. We now prove the claim for $h=h'g_{n-1}g_{n-2}...g_{n-s}$ by induction on $s$, with the
case for $s=1$ just shown. Using the relation $(g_{n-s+1}g_{n-s}-g_{n-s+1})e_{(n)}=(g_{n-s}g_{n-s+1}-g_{n-s})e_{(n)}$,
we obtain
$$e_{(n)}h'g_{n-1}g_{n-2}...g_{n-s} e_{(n)}= $$
$$=e_{(n)}(h'g_{n-s}-h')g_{n-1}g_{n-2}...g_{n-s+1}e_{(n)} 
+ e_{(n)}(h'g_{n-s})g_{n-1}g_{n-2}...g_{n-s+2}e_{(n)}.$$
The claim now holds for each summand on the right hand side by induction assumption. After applying $\Theta_n$ to it,
it can be easily shown that it is equal to $e_{(n)}\Theta(h)e_{(n)}$.

\subsection{Trace property of $\phi$} We will use the following simple observation. Let $A$ be a semisimple algebra,
and let $I\subset A$ be a two-sided ideal. If $\psi: A\to \C$ is a functional satisfying $\psi(cd)=\psi(dc)$ for all $c,d\in I$,
then we also have $\psi(ca)=\psi(ac)$ for all $a\in A$ and $c\in I$. Indeed, we can write $a=a_I+a_J$ with $a_I\in I$
and $a_J\in J$ where $J$ is a two-sided ideal of $A$ such that $IJ=0$. The claim follows from $a_Jc=0=ca_J$.

\begin{theorem}\label{Markovtheorem}
The functional $\phi$ satisfies the Markov property  $\phi(cg_{n-1})=\phi(c)\phi(g_{n-1})$
for all $c\in C_{n-1}$ and the
trace property $\phi(cd)=\phi(dc)$ for all $c,d\in C_n$. Hence there exists a trace on $C_n$ satisfying
Condition \ref{Markovprop}.
\end{theorem}

$Proof.$ The first claim follows from Lemma \ref{Markovlemma},(b).
We will prove the second claim by induction on $n$, which is certainly true for the abelian algebra $C_1$.
 Let $I_r=\bigoplus_{s\geq r} H_ne_{(s)}H_n$. We will prove that the restriction of $\phi$ to $I_r$ satisfies
the trace property by downwards induction. 
We define the functional $\al: H_n\to \C$ by 
$$e_{(n)}he_{(n)}=\al(h)e_{(n)}.$$
It follows from Lemma \ref{Markovlemm1} and Lemma \ref{Markovlemma2}  that
$${[N]^n}\phi(h_1e_{(n)}h_2)=Tr(\Phi(\Theta_n(h_1)E^{\otimes n}\Phi(\Theta_n(h_2))=$$
$$= Tr(E^{\otimes n}\Phi(\Theta_n(h_2h_1))(E^{\otimes n})= Tr(E^{\otimes n}\Phi(h_2h_1)E^{\otimes n})= \al(h_2h_1).$$
 Then we calculate
$$[N]^n\phi(a_1e_{(n)}a_2b_1e_{(n)}b_2)=
[N]^n\al(a_2b_1)\phi(a_1e_{(n)}b_2) =\al(a_2b_1)\al(b_2a_1).$$
One calculates in the same way that also $[N]^n\phi(b_1e_{(n)}b_2a_1e_{(n)}a_2)=\al(a_2b_1)\al(b_2a_1)$.
This proves the claim for $r=n$.
For the induction step, we first prove 
\begin{equation}\label{Markovind}
\phi(he_{(s)})=\phi(e_{(s)}he_{(s)})\quad {\rm for\ all\ }h\in H_n, s\geq r.
\end{equation}
For $s>r$, this is clear as the restriction of $\phi$ to $H_ne_{(s)}H_n$ is a trace by induction assumption.
For $s=r$, we will prove the claim by induction on $m\geq r$, which has already been proved
for $m=r$. For the induction step, it suffices to prove the claim for elements of the form $h=h_1g_mh_2$, $h_1,h_2\in H_m$.
But by Markov property, we have
$$\phi(h_1g_mh_2e_{(r)})=\phi(g_m)\phi(h_1h_2e_{(r)})=\phi(g_m)\phi(e_{(r)}h_1h_2e_{(r)})=\phi(e_{(r)}h_1g_mh_2e_{(r)}),$$
where we used the induction assumption for $h_1h_2\in H_m$.

Let now $h\in H_n$ and $k\in I_{r+1}$ such that $ke_{(r)}=e_{(r)}k=k$. Then it follows from \ref{Markovind} that
\begin{equation}\label{traceproof}
\phi(hk)=\phi(hke_{(r)})=\phi(e_{(r)}hke_{(r)})=\phi(e_{(r)}he_{(r)}k),
\end{equation}
where we used that $\phi$ is a trace on $I_{r+1}$.
Let now $a_i,b_i\in H_n$, $ i=1,2$. Then we can write $e_{(r)}a_2b_1e_{(r)}=(h'+k)e_{(r)}=e_{(r)}(h'+k)$
with $h'\in H_{r+1,n}$ and $k\in I_{r+1}$.
We then calculate, using \ref{traceproof}
$$\phi(a_1e_{(r)}a_2b_1e_{(r)}b_2)=\phi(a_1e_{(r)}(h'+k)e_{(r)}b_2)=$$
$$=\phi(e_{(r)}b_2a_1e_{(r)}(h'+k))
=\phi(e_{(r)}b_2a_1e_{(r)}a_2b_1e_{(r)}).
$$
We similarly calculate 
$$\phi(b_1e_{(r)}b_2a_1e_{(r)}a_2)=\phi(e_{(r)}a_2b_1e_{(r)}b_2a_1e_{(r)}).$$
The trace property now follows from the induction assumption, using the fact that $e_{(r)}C_ne_{(r)}\cong C_{n-r}$,
see Proposition \ref{basicstructure}

\section{Conclusions and future research}\label{conclusions}

\subsection{Historical context} 
When the author of this paper visited Columbia University as a postdoc in the first half of 1986, Vaughan suggested
as a project with my host, J. Birman, that we try to find an algebraic interpretation of the Kauffman polynomial.
This resulted in the definition of a new algebra (independently discovered by J. Murakami) which turned out to be
a $q$-deformation of Brauer's centralizer algebra, see \cite{Brauer}, \cite{BBMW}. This algebra was subsequently
used to construct subfactors of type $BCD$, among other applications. A different $q$-deformation of Brauer's centralizer
algebra was found in \cite{WqBrauer}, see also \cite{Mo}. It contained the Hecke algebra $H_n$ as a subalgebra.
It was shown in \cite{Wfusion} that it could also be used to construct subfactors, as well as module categories of $Rep(U_q\sl_N)$ and of the
related fusion tensor categories $\overline{Rep}(U_q\sl_N)$ for $q$ a root of unity (they are often also referred to
as $SU(N)_k$). This will be sketched below. In particular, we could explicitly calculate
the indices and first principal graphs of these subfactors. This, in turn, also allows us to give an explicit description of the algebras
corresponding to these module categories. So while the current paper is purely algebraic, it is closely related to 
research in which Vaughan was interested. In particular, the idea of a Markov trace which plays a crucial role
in finding the relations for the algebras $C_n$ goes back to him. It should also be noted that the first module categories
for fusion categories related to $SU(2)$ were already constructed by Vaughan and his collaborators in \cite{GHJ}.

\subsection{Markov traces, module categories and subfactors} 
It was shown in \cite{WHecke} that the quotient $\overline{H}_n(q)$ of $H_n(q)$
modulo the annihilator of the Markov trace $tr$ is semisimple for all $n$, with $tr$ as in Remark
\ref{MarkovHecke}. We expect the same to be true for the quotient
$\overline{C}_n(q)$ modulo the annihilator ideal of its extension, which was shown to exist in Theorem \ref{Markovtheorem}.
 E.g. it is not hard to see that for $q$ not a root of unity the quotient is isomorphic to
the image of $C_n(q)$ in its representation in $\End(V^{\otimes n})$, see Section \ref{tensorsection}. Assuming this,
the construction of the module category goes as follows:

It was shown in \cite{KW} (see also \cite{TW} for a variation of this construction) that $Rep(U_q\sl_N)$ for $q$ not
a root of unity and $\overline{Rep}(U_q\sl_N)$ for $q$ a root of unity can be reconstructed from the quotients
$\overline{H}_n(q)$  modulo the annihilator ideal of a suitable version of the Markov trace. 
Here the objects are given by idempotents of $\overline{H}_n(q)$. We similarly define
the module category $\M$ whose objects are  idempotents in $\overline{C}_n(q)$. If 
$p_M\in \overline{C}_m(q)$ and $p_H\in \overline{H}_n(q)$ are  idempotents, we define the module action by
$$p_M\otimes p_H\ :=\ p_Msh_m(p_H),$$
where the algebra homomorphism  $sh_m: H_n(q)\to C_{n+m}(q)$ is defined by $sh_m(g_i)=g_{i+m}\in C_{n+m}(q)$. It follows from the relations
that $p_M\otimes p_H$ is an idempotent in $\overline{C}_{n+m}(q)$. 
Finally, if the quotients $\overline{H}_n(q)$ and $\overline{C}_n(q)$ allow compatible $C^*$ structures, we can construct
subfactors $\Na\subset \M$ from the inclusions $\lim_{n\to\infty} \overline{H}_n(q)\subset\overline{C}_n(q)$
following the procedure in \cite{WHecke}, Section 1.

\subsection{Results in \cite{WqBrauer} and \cite{Wfusion}}  We give an outline of the results in these papers
which give a good idea of the results to be expected in the approach outlined in the previous subsection.
A $q$-version $Br_n(q)$ of Brauer's centralizer algebra (see \cite{Brauer}) was defined in \cite{WqBrauer} by again adding one more generator $e$ to
the generators of the Hecke algebras $H_n(q)$. As in this paper, the relations were forced by the condition that
the extension of the
Markov trace on $H_n(q)$ to the algebras $Br_n(q)$ satisfy an analog of the Markov condition \ref{Markovprop}.
Subfactors were constructed from these algebras as outlined in the previous section. Their indices and first principal
graphs are given in \cite{Wfusion} Sections 3F and 3G.  Instead of copying the results there, we just state an easy
consequence of these results which only appears implicitly in \cite{Wfusion}:

Let $q=e^{\pi i/(N+k)}$. Then the category constructed from the quotients $\overline{H}_n(q)$ is equivalent
to the fusion category $SU(N)_k$ (or $\overline{Rep}(U_q\sl_N)$ in the notation of this paper). For simplicity,
we assume a trivial twist (see \cite{KW} or \cite{TW} for details). It is well-known that the simple objects of $SU(N)_k$ are labeled
by the Young diagrams $\la$ with $\leq N-1$ rows such that $\la_1\leq k$. Recall that a module category
over a tensor category can be defined via an algebra object in the given tensor category (see \cite{ostrik}). 
We will reformulate the following theorem in a somewhat more conceptual way in 
 Remark \ref{concremark}.

\begin{theorem}\label{exampletheorem} Let $N$ be even. Then $SU(N)_k$ has an algebra object $A=Ind_{Ad}(\1)$,
where $Ind_{Ad}(\1)$ is the direct sum of simple objects $V_\la$ such that $N|\ |\la|$ and the number of boxes in each column of $\la$ is
even.
\end{theorem}

$Proof.$ We consider the inclusion of von Neumann factors $\Na\subset \M$ constructed in  \cite{Wfusion}, 
Theorem 3.4  for case (c) listed before that theorem. It follows from the explicit description of its
principal graph  in \cite{Wfusion}, 
Section 3G that the von Neumann algebra $\M$, viewed as an $\Na-\Na$ bimodule, decomposes into a direct sum of
simple  $\Na-\Na$  bimodules labeled by exactly the Young diagrams which appear in $Ind_{Ad}(\1)$.
As $\M$ has a multiplication, it follows that $A$ in the statement is an algebra object in the category of
 $\Na-\Na$ bimodules. It is known that this category is equivalent to $Ad(SU(N)_k$, the subcategory of $SU(N)_k$ 
whose simple objects are labeled by Young diagrams $\la$ such that $N | \ |\la|$ (see e.g. \cite{WCstar}, Theorem 4.4).

\begin{remark}\label{concremark} 1. It is well-known that the restriction to $Sp(N)$ of a simple $SU(N)$-module labeled by the Young diagram $\la$ 
contains the trivial representation of $Sp(N)$ if and only if the number of boxes in each column of $\la$  is even. Hence 
the algebra in Theorem \ref{exampletheorem} can be viewed 
 as a natural analog in $SU(N)_k$ of the induction of the trivial representation of $PSp(N)$ to $PSU(N)$.
It would seem plausible that similar algebras exist which correspond to inducing the trivial representation of $PSp(N)$ or $Sp(N)$ 
to quotients of $SU(N)$ modulo a subgroup of its center. This seems to be compatible with Edie-Michell's classification
results of module categories of $SU(N)_k$, see \cite{Cain} and \cite{Caincomm}. This possible generalization of our result
became evident after conversations with Edie-Michell.

2. Similarly, one can construct module categories and algebras from cases (a) and (b) before \cite{Wfusion}, Theorem 3.4. 
As they are related to embeddings of the full orthogonal group $O(N)$, we would need as larger group the group $SU(N)\times \Z/2$
of unitary matrices $u$ with $| \det(u)| = 1$. It should be possible to obtain module categories of
$SU(N)_k$ from this via a $\Z/2$ orbifold construction. One obtains algebra objects for these cases in the same way as it was done in
Theorem \ref{exampletheorem}. Again, one would expect algebra objects and module categories corresponding to each quotient
group $SU(N)/Z$, where $Z$ is a subgroup of the center of $SU(N)$.

3. A complete realization of all module categories for all fusion tensor categories of type $SU(3)_k$
has been given in \cite{EP}. Using their results, one can show that the general approach outlined here also
works  in the setting of this paper for 
the special case $N=3$,
i.e. for the embedding of $Sp(2)\subset SL(3)$. Indeed, the explicit calculations in \cite{EP} were useful in the initial phase of
finding relations for our algebras.  Here the algebra object coming from the subfactor
constructed there would be the direct sum of all simple objects in $SU(3)_k$
labeled by Young diagrams $\la$ with $3|\ |\la|$. These subfactors seem to be closely connected to subfactors constructed by F. Xu in \cite{Xu}.
\end{remark}
		
\subsection{Classification of module categories of WZW-fusion categories}  A lot of progress in classifying module categories
has recently been made by Edie-Michell \cite{Cain}, building on the works of Ocneanu, Gannon, Schopieray, Evans and Pugh,  and others.
Very roughly speaking these module categories can be divided into exceptional and non-exceptional
module categories. It appears that we can find realizations for all non-exceptional module categories of
fusion tensor categories of type $SU(N)_k$ using the construction sketched in this paper and its generalizations in
Remark \ref{concremark} together with the orbifold construction. This was done in collaboration with Edie-Michell.
  It would also be interesting to find out whether non-exceptional
module categories of fusion categories of other Lie types could similarly be realized via the constructions in this paper
for certain subgroups in connection with orbifolds.

\subsection{Co-ideal subalgebras} As mentioned in the introduction, module categories of a Drin- feld-Jimbo quantum group $U_q\g$ can be defined
for the sub-Lie algebra $\h$ consisting of the fixed points of an order 2 Lie algebra automorphism, see \cite{Le} and \cite{NS}. 
It would be interesting to see whether our (proposed) module category of $U_q\sl_N$ could be realized by a suitable co-ideal deformation
of the universal enveloping algebras $U\sp_{N-1}\subset U\sl_N$.

\end{document}